\documentclass[11pt]{amsart}
\usepackage[T1]{fontenc}
\usepackage{lmodern}
\usepackage[margin=25mm,headheight=14pt]{geometry}
\usepackage{amsaddr} % Moves affiliations/addresses to the first page
\usepackage{amsmath,amssymb,amsthm,mathtools,mathrsfs}
\usepackage{booktabs,array,enumitem,fancyhdr,microtype}
\usepackage{xcolor}
\usepackage[unicode,colorlinks=true,linkcolor=blue!40!black,citecolor=blue!40!black,urlcolor=blue!40!black]{hyperref}
\hypersetup{pdftitle={A counterexample to logarithmic choosability for bipartite graphs},pdfsubject={A proof by weighted counting and a fixed-point argument}}
\allowdisplaybreaks[2]
\numberwithin{equation}{section}
\newtheorem{theorem}{Theorem}[section]
\newtheorem{lemma}[theorem]{Lemma}
\newtheorem{proposition}[theorem]{Proposition}
\newtheorem{corollary}[theorem]{Corollary}
\theoremstyle{definition}
\newtheorem{definition}[theorem]{Definition}
\newtheorem{conjecture}[theorem]{Conjecture}
\newtheorem{example}[theorem]{Example}

\theoremstyle{remark}
\newtheorem{remark}[theorem]{Remark}
\newcommand{\ch}{\operatorname{ch}}
\newcommand{\calF}{\mathcal F}
\newcommand{\calG}{\mathcal G}
\newcommand{\calH}{\mathcal H}
\newcommand{\calS}{\mathcal S}
\newcommand{\R}{\mathbb R}
\newcommand{\avg}[1]{\langle #1\rangle}
\setlist[enumerate]{itemsep=0.3em,topsep=0.4em}
\title{\bfseries A Counterexample to the Alon-Krivelevich Conjecture and a New Lower Bound for Choosability of Bipartite Graphs}
\author{Junran Yu$^{1}$ and Gexin Yu$^{2}$}

\address{$^{1}$\small School of Mathematics, Nanjing Normal University, Nanjing, Jiangsu, China.\\
$^{2}$\small Department of Mathematics, William \& Mary, Williamsburg, VA 23185, USA. 
}

\thanks{Email: gyu@wm.edu}

\begin{document}
\maketitle

\begin{abstract}
For an integer $k\ge 2$ and $r^2>2(2k-1)\ln (2k-1)$, we construct a bipartite graph $G$ whose two vertex classes are disjoint copies of the family of all $k$-subsets of a $(2k-1)$-element set, two vertices being adjacent when their sets intersect in at most $r$ elements.
Assigning each vertex its corresponding subset as its  color list, $G$ is not colorable, and thus its choice number $\ch(G)\ge k+1$. This implies that there are graphs with maximum degree $\Delta$ and
$$\ch(G)\ge\Bigl(\frac1{32}-o(1)\Bigr)\frac{(\ln\Delta)^2}{(\ln\ln\Delta)^3}.$$ When $k\to\infty$,
consequently $\ch(G)/\ln\Delta$ tends to infinity, which disproves the Alon--Krivelevich conjecture. In order to show the uncolorability, we use tools from probability (conditional product measure, nested
coupling, and Doob martingale), which we write in an elementary language of weighted counting, and topology (Brouwer fixed-point theorem). 
\end{abstract}

\medskip
\noindent\textbf{Keywords.} Bipartite graph; list coloring;  weighted counting; Brouwer fixed-point theorem.

\section{Introduction}\label{sec:intro}

\subsection{List coloring and the Alon--Krivelevich conjecture}
All graphs in this paper are finite and simple.
A \emph{list assignment} $L$ on a graph $G=(V,E)$ assigns a finite set $L(v)$ of available colors to every vertex $v\in V$.
A \emph{proper $L$-coloring} is a function $\varphi$ satisfying
$\varphi(v)\in L(v)$ with $v\in V$ and
$\varphi(u)\ne \varphi(v)$ with $uv\in E$.
The graph is \emph{$q$-choosable} if every assignment of lists of size at least $q$ has a proper list coloring. Its \emph{choosability}, also called its \emph{list chromatic number}, is
$\ch(G)=\min\{q\in\mathbb Z_{\ge1}:G\text{ is }q\text{-choosable}\}$.

%For a nonempty graph this minimum exists: lists of size $|V|$ permit a greedy coloring in any vertex order, since at each step at most $|V|-1$ previously colored vertices can forbid colors. The graph with no vertices is $1$-choosable by the empty coloring. Choosability is monotone in the list size. Indeed, if $G$ is $q$-choosable and every vertex receives a list of size $q+1$, choose a $q$-element sublist at each vertex and apply $q$-choosability to the resulting assignment. Thus an uncolorable assignment of lists of size $k$ implies $\ch(G)\ge k+1$.

Although the ordinary chromatic number of a bipartite graph is at most two, its choosability can be arbitrarily large.
Alon and Krivelevich~\cite[Conjecture 5.1]{AK1998} made the following conjecture.

\begin{conjecture}[Alon and Krivelevich~\cite{AK1998}]
There exists an absolute constant $C>0$ such that
\begin{equation}\label{eq:conjecture}
\ch(G)\le C\ln\Delta(G)
\end{equation}
for every finite simple bipartite graph $G$ with maximum degree $\Delta(G)\ge2$.
\end{conjecture}

Bradshaw, Mohar, and Stacho~\cite[Theorem 1.2]{BMS2024} proved that every bipartite graph of sufficiently large maximum degree satisfies
$\ch(G)<0.797\Delta(G)/\ln\Delta(G)$. The authors and Ninh~\cite{YNY} extended the Bradshaw--Mohar--Stacho framework to list coloring of near-bipartite graphs and obtained the bound $0.77\Delta(G)/\ln\Delta(G)$.
The results below give an explicit family for which \eqref{eq:conjecture} fails.

\subsection{Notation}
For a positive integer $N$, write $[N]=\{1,\ldots,N\}$, and put $[0]=\varnothing$.
For a finite set $\Omega$ and an integer $j$, let
$\binom\Omega j=\{U\subseteq \Omega:|U|=j\}$,
which we call the \emph{$j$-th layer} of $\Omega$, it is empty if $j<0$ or $j>|\Omega|$.
The notation $\binom Nj$ denotes the usual numerical binomial coefficient.
Throughout, $\ln$ is the natural logarithm and $e$ is its base.

\subsection{The construction and the results}
\begin{definition}\label{def:graph}
Let $k\ge2$ and $r\ge0$ be integers, and put $P=[2k-1]$. The graph $G(k,r)$ has the disjoint vertex classes
$A=\{a_S:S\in\tbinom Pk\}$ and  $B=\{b_T:T\in\tbinom Pk\}$
and the edge of $G(k,r)$ is defined as follows:
\begin{equation}\label{eq:edge}
a_Sb_T\in E(G(k,r))\Longleftrightarrow |S\cap T|\le r.
\end{equation}
Its \emph{standard lists} are
\begin{equation}\label{eq:lists}
L(a_S)=S\ \text{and} \ L(b_T)=T.
\end{equation}
\end{definition}

Two $k$-subsets of a $(2k-1)$-element set always intersect, and a typical pair intersects in about $k/2$ elements. By~\eqref{eq:edge}, $a_S$ is adjacent only to the vertices $b_T$ whose set $T$ is nearly the complement of $S$. When $r$ is much smaller than $k$ such pairs are rare, which keeps the degree small, while the lists have size $k$.

\begin{theorem}[Uncolorability criterion]\label{thm:criterion}
Let $k\ge2$ and $r\ge0$ be integers, and put $n=2k-1$. If
\begin{equation}\label{eq:threshold}
r^2>2n\ln n,
\end{equation}
then $G(k,r)$ has no proper coloring from its standard lists.
\end{theorem}

\begin{corollary}\label{thm:main}
For an integer $k\ge2$, let $n=2k-1$,
$r=\bigl\lfloor\sqrt{2n\ln n}\bigr\rfloor+1$,
$G=G(k,r)$ and
$\Delta=\Delta(G)$.

\begin{enumerate}[label=\textup{(\alph*)}]
\item The standard lists of $G$ admit no proper coloring, hence $\ch(G)\ge k+1$.
\item As $k\to\infty$,
\begin{equation}\label{eq:optimal}
\ch(G)\ge\Bigl(\frac1{32}-o(1)\Bigr)\frac{(\ln\Delta)^2}{(\ln\ln\Delta)^3}.
\end{equation}
In particular $\ch(G)/\ln\Delta\to\infty$, so no absolute constant $C$ makes~\eqref{eq:conjecture} valid for all finite simple bipartite graphs.
\item The order in~\eqref{eq:optimal} is the limit of Theorem~\ref{thm:criterion} for these graphs: if $r'$ is any integer with $(r')^2>2n\ln n$ and $\Delta'=\Delta(G(k,r'))$, then $k\le(1+o(1))(\ln\Delta')^2/(\ln\ln\Delta')^3$ as $k\to\infty$, uniformly in the choice of $r'$.
\end{enumerate}
\end{corollary}

Part~(a) is Theorem~\ref{thm:criterion} applied to this $r$, since $r>\sqrt{2n\ln n}$. For fixed $k$, the bound $\ch\ge k+1$ does not depend on $r$, while the edge sets of $G(k,r)$, and hence $\Delta(G(k,r))$, are nondecreasing in $r$; so this $r$, the least integer satisfying~\eqref{eq:threshold}, is the best choice for this method. Part~(c) concerns only the lower bound produced by Theorem~\ref{thm:criterion}; it says nothing about the true choosability of $G(k,r)$, which could be larger, nor about the true least $r$ for which $G(k,r)$ is uncolorable from its standard lists.

\subsection{Organization and an outline of the proof}
Section~\ref{sec:construction} establishes the properties of $G(k,r)$ other than uncolorability: it is a finite simple bipartite regular graph, and its degree is computed exactly and bounded. It then derives Corollary~\ref{thm:main} from Theorem~\ref{thm:criterion}.

Sections~\ref{sec:coloring}--\ref{sec:separation} prove Theorem~\ref{thm:criterion} by contradiction. We assume that a proper coloring $\varphi$ exists and follow its consequences.
\begin{enumerate}[label=\textbf{Step \arabic*.},leftmargin=*]
\item (Section~\ref{sec:coloring}) The coloring defines two selection rules. They partition the $k$-sets and the $(k-1)$-sets into color classes $\calF_i$ and $\calG_i$. Properness says that, for every color $i$, every set of $\calF_i$ is far from every set of $\calG_i$, and hence far from every $k$-set containing a member of $\calG_i$. Here ``far'' means differing in at least $r$ elements.
\item (Section~\ref{sec:weights}) To turn ``far apart'' into a contradiction, we must know that some color has \emph{large} classes on both sides. We measure families of sets by weighted proportions: a vector $x\in[0,1]^n$ gives the set $U$ the weight $\prod_{i\in U}x_i\prod_{i\notin U}(1-x_i)$. Plain proportions, the case $x=(\tfrac12,\ldots,\tfrac12)$, are not enough, because a coloring can make the shared colors rare.
\item (Section~\ref{sec:balance}) We choose the weights \emph{from the coloring}: a fixed-point argument yields $x^*$ such that every color has the same weighted share among $k$-sets as among $(k-1)$-sets. Some color $c$ then has share at least $1/n$ on both sides. 
\item (Section~\ref{sec:lift}) The $(k-1)$-sets of color $c$ are lifted to the $k$-sets containing them, without losing weighted share.
\item (Section~\ref{sec:separation}) Two families of $k$-sets with weighted share at least $1/n$ always contain two sets that differ in at most $\sqrt{2n\ln n}$ elements.
\end{enumerate}

\begin{proof}[Proof of Theorem~\ref{thm:criterion}]
Under the standing assumption of Section~\ref{sec:coloring}, let $x^*$ and $c$ be as in Lemma~\ref{lem:balance}, and consider the two families
$\calF=\calF_c$ and $\calH=\calH_c=\nabla\calG_c$
in the $k$-th layer of $P$. We collect what has been proved about them.
\begin{enumerate}[label=\textup{(\roman*)}]
\item $\nu_{x^*,k}(\calF)\ge1/n$, by Lemma~\ref{lem:balance}.
\item $\nu_{x^*,k}(\calH)\ge1/n$, by Corollary~\ref{cor:lift}.
\item $d_J(S,U)\ge r$ for all $S\in\calF$ and $U\in\calH$, by Lemma~\ref{lem:far}.
\end{enumerate}
Lemma~\ref{lem:separation} with $N=n$, $j=k$, $x=x^*$, and $\rho=1/n$ now gives
$r\le\sqrt{2n\ln n}$.
Squaring is valid because both sides are nonnegative, and the result contradicts~\eqref{eq:threshold}. Hence no proper coloring exists.
\end{proof}

The coloring $\varphi$ enters the argument only through the color classes in Step~1 and the choice of $x^*$ in Step~3. Everything proved in Sections~\ref{sec:weights}, \ref{sec:lift} and~\ref{sec:separation} is a general fact about weights on the subsets of a finite set. The only external ingredient is Brouwer's fixed-point theorem, used in Step~3.

\section{The construction and its properties}\label{sec:construction}
The graph $G(k,r)$ is bipartite with $2\binom{2k-1}k$ vertices; since $a_S$ and $b_T$ are distinct vertices even when $S=T$, it has no loops and no multiple edges. All standard lists have size $k$.

Fix $S\in\binom Pk$. As $|P\setminus S|=k-1$, every $k$-set $T$ meets $S$, and a $k$-set $T$ with $|S\cap T|=t$ is determined by the independent choices of $S\cap T\subseteq S$ and $T\setminus S\subseteq P\setminus S$. The same count applies to the vertices $b_T$, so $G(k,r)$ is regular of degree
\begin{equation}\label{eq:exact-degree}
\Delta(G(k,r))=\sum_{t=1}^{\min(r,k)}\binom kt\binom{k-1}{k-t}=\sum_{t=1}^{\min(r,k)}\binom kt\binom{k-1}{t-1}.
\end{equation}
The edge sets grow with $r$, so $\Delta(G(k,r))$ is nondecreasing in $r$. Let $1\le r\le k/2$. For $1\le t\le r$ we have $\binom{k-1}{t-1}=\frac tk\binom kt\le\binom kt\le\binom kr$, the last inequality because $\binom k{s+1}/\binom ks=\frac{k-s}{s+1}\ge1$ whenever $2s+1\le k$. Stirling's formula gives $r!\ge(r/e)^r$, hence $\binom kr\le k^r/r!\le(ek/r)^r$; and $\binom kr=\prod_{i=0}^{r-1}\frac{k-i}{r-i}\ge(k/r)^r$. Since the term $t=r$ of~\eqref{eq:exact-degree} is at least $\binom kr$,
\begin{equation}\label{eq:degree-bounds}
\Bigl(\frac kr\Bigr)^r\le\Delta(G(k,r))\le r\binom kr^2\le r\Bigl(\frac{ek}r\Bigr)^{2r}\qquad(1\le r\le k/2).
\end{equation}

\begin{proof}[Proof of Corollary~\ref{thm:main} from Theorem~\ref{thm:criterion}]
All $o(1)$ terms refer to $k\to\infty$. Write $L=\ln\Delta$ and $L'=\ln\Delta'$. Since $2n\ln n=(4+o(1))k\ln k$ and $\sqrt{2n\ln n}<r\le\sqrt{2n\ln n}+1$, we have $r=(2+o(1))\sqrt{k\ln k}$; hence, for large $k$, $1\le r\le k/2$, and $\ln(k/r)$, $\ln(ek/r)$ both equal $(\frac12+o(1))\ln k$. By~\eqref{eq:degree-bounds},
$$
(1+o(1))\sqrt k(\ln k)^{3/2}\le r\ln\frac kr\le L\le\ln r+2r\ln\frac{ek}r=(2+o(1))\sqrt k(\ln k)^{3/2}.$$
Hence $\ln L=(\frac12+o(1))\ln k$. Part~(a) gives $\ch(G)\ge k+1$, and the upper bound on $L$ gives $k\ge L^2/((4+o(1))(\ln k)^3)=(\frac1{32}-o(1))L^2/(\ln L)^3$. This proves~(b), and $\ch(G)/L\to\infty$ because $L\to\infty$.

For~(c), the integer $r'$ exceeds $\sqrt{2n\ln n}$, so $r'\ge r$. By monotonicity in $r$ and~\eqref{eq:degree-bounds}, $L'\ge r\ln(k/r)=(1+o(1))\sqrt k(\ln k)^{3/2}$, where the $o(1)$ depends only on $k$. Also $\Delta'\le\binom{2k-1}k\le4^k$, so $\ln L'\le\ln k+1$. Therefore $k\le(1+o(1))(L')^2/(\ln k)^3\le(1+o(1))(L')^2/(\ln L')^3$.
\end{proof}

It remains to prove Theorem~\ref{thm:criterion}, this occupies Sections~\ref{sec:coloring}--\ref{sec:separation}.

\section{From a proper coloring to far-apart color classes}\label{sec:coloring}

\subsection{The standing assumption}
Throughout Sections~\ref{sec:coloring}--\ref{sec:separation}, fix integers $k\ge2$ and $r\ge0$ satisfying~\eqref{eq:threshold}, let $n=2k-1$ and $P=[n]$, and suppose for contradiction that
$\varphi$ is a proper coloring of $G(k,r)$ from its standard lists.
Sections~\ref{sec:weights}, \ref{sec:lift} and~\ref{sec:separation} also contain general lemmas that do not use $\varphi$, these are stated for an arbitrary finite ground set.

\subsection{Two selection rules}
Define
\begin{equation}\label{eq:rules}
f(S)=\varphi(a_S)\quad\Bigl(S\in\tbinom Pk\Bigr)
\ \ \text{and } \
g(R)=\varphi(b_{P\setminus R})\quad\Bigl(R\in\tbinom P{k-1}\Bigr).
\end{equation}
The second definition is legitimate since $|P\setminus R|=(2k-1)-(k-1)=k$, so $b_{P\setminus R}$ is a vertex of $B$.
The list restrictions give
$f(S)\in S$, $g(R)\in P\setminus R$ that is $g(R)\notin R$.
Thus $f$ selects an element \emph{inside} each $k$-set, and $g$ selects an element \emph{outside} each $(k-1)$-set.

\subsection{Color classes}
 Let
$\calF_i=\{S\in\tbinom Pk:f(S)=i\}$ and
$\calG_i=\{R\in\tbinom P{k-1}:g(R)=i\}$ for each color $i\in P$.
The families $\calF_1,\ldots,\calF_n$ partition the $k$-th layer of $P$, and $\calG_1,\ldots,\calG_n$ partition the $(k-1)$-th layer. Since $f(S)=i$ for $S\in\calF_i$ and $g(R)=i$ for $R\in\calG_i$, the list restrictions give
\begin{equation}\label{eq:class-containment}
i\in S\ \ \text{for every } S\in\calF_i \ \ \text{and } \
i\notin R\ \ \text{for every } R\in\calG_i.
\end{equation}

\subsection{Properness as a distance condition}
For $U,W\in\binom\Omega j$, where $\Omega$ is a finite set, define the \emph{Johnson distance}
\begin{equation}\label{eq:johnson}
d_J(U,W)=|U\setminus W|=|W\setminus U|
=\tfrac12|U\mathbin\triangle W|,
\end{equation}
the number of elements that must be exchanged to turn $U$ into $W$.
The first two quantities agree because $|U|=|W|$. Nonnegativity and symmetry are immediate, and $d_J(U,W)=0$ holds exactly when $U=W$. For the triangle inequality, the inclusion
$U\setminus W\subseteq(U\setminus V)\cup(V\setminus W)$
gives $d_J(U,W)\le d_J(U,V)+d_J(V,W)$. Thus $d_J$ is a metric on each layer, and $d_J(U,W)\le j$.

\begin{lemma}[Same color forces distance]\label{lem:distance}
If $S\in\calF_i$ and $R\in\calG_i$ for the same color $i$, then $|S\setminus R|\ge r+1$.
\end{lemma}
\begin{proof}
The vertices $a_S$ and $b_{P\setminus R}$ both receive color $i$.
Since $\varphi$ is proper, they are not adjacent, so~\eqref{eq:edge} gives $|S\cap(P\setminus R)|>r$.
Both sides are integers, and $S\cap(P\setminus R)=S\setminus R$.
\end{proof}

The two classes of a color live on different layers. To compare them on one layer, we enlarge the $(k-1)$-sets. For $\calG\subseteq\binom\Omega{j-1}$, let
$
\nabla\calG=\Bigl\{U\in\tbinom\Omega j:\ U\supseteq R\text{ for some }R\in\calG\Bigr\}$
be the \emph{upper shadow} of $\calG$, and let
$\calH_i=\nabla\calG_i\subseteq\tbinom Pk$ with $i\in P$.

\begin{lemma}[Far-apart classes]\label{lem:far}
For every color $i$, every $S\in\calF_i$, and every $U\in\calH_i$, we have $d_J(S,U)\ge r$.
\end{lemma}
\begin{proof}
Choose $R\in\calG_i$ with $R\subseteq U$. By Lemma~\ref{lem:distance}, $|S\setminus R|\ge r+1$.
Also $|U\setminus R|=1$. Because $R\subseteq U$, every element of $S\setminus R$ either lies outside $U$ or belongs to $U\setminus R$, and these two possibilities are disjoint:
$S\setminus R=(S\setminus U)\sqcup\bigl(S\cap(U\setminus R)\bigr)$.
Taking cardinalities and using $|S\cap(U\setminus R)|\le1$,
$d_J(S,U)=|S\setminus U|=|S\setminus R|-|S\cap(U\setminus R)|\ge(r+1)-1=r$.\qedhere
\end{proof}

\subsection{What remains to be shown}
Lemma~\ref{lem:far} is the only consequence of properness used below. To contradict it, it suffices to find one color $c$ for which $\calF_c$ and $\calH_c$ are both ``large''; Section~\ref{sec:separation} shows that two large families of $k$-sets contain two sets at Johnson distance at most $\sqrt{2n\ln n}<r$.

\begin{remark}[A common color always exists]\label{rem:common}
Some color is used by both rules, by a counting argument. Let $C_f$ and $C_g$ be the sets of colors used by $f$ and by $g$. Every $k$-set $S$ contains $f(S)\in C_f$, so $C_f$ meets every $k$-set; hence $|P\setminus C_f|\le k-1$, that is, $|C_f|\ge k$. If $|C_g|\le k-1$, choose a $(k-1)$-set $R\supseteq C_g$; then $g(R)\in C_g\subseteq R$, which is impossible. Hence $|C_g|\ge k$. Since $|C_f|+|C_g|\ge2k>|P|$, the two sets share a color. This argument gives no information about the \emph{sizes} of the corresponding classes, which is what the rest of the proof provides.
\end{remark}

\section{Measuring families of sets by weights}\label{sec:weights}
This section is independent of $\varphi$. Throughout it, $\Omega$ is a finite set and $x=(x_i)_{i\in\Omega}\in[0,1]^\Omega$.

\subsection{Product weights}
\begin{definition}[Product weights]\label{def:weights}
The \emph{$x$-weight} of a set $U\subseteq\Omega$ is
$$\pi_x(U)=\prod_{i\in U}x_i\prod_{i\in\Omega\setminus U}(1-x_i),$$
an empty product being $1$. For a family $\calF$ of subsets of $\Omega$, put $\pi_x(\calF)=\sum_{U\in\calF}\pi_x(U)$.
\end{definition}

The parameter $x_i$ is a degree of preference for element $i$: the closer $x_i$ is to $1$, the heavier the sets containing $i$ are, relative to the sets avoiding it. If $x_i=1$, every set avoiding $i$ has weight zero, so $i$ is \emph{forced}; if $x_i=0$, every set containing $i$ has weight zero, so $i$ is \emph{forbidden}. If every $x_i$ equals $1/2$, every set has weight $2^{-|\Omega|}$.

\begin{lemma}[Basic properties of the weights]\label{lem:weights}
Let $J_0=\{i\in\Omega:x_i=0\}$, $J_1=\{i\in\Omega:x_i=1\}$, and $Q=\Omega\setminus(J_0\cup J_1)$.
\begin{enumerate}[label=\textup{(\alph*)}]
\item $\pi_x(U)\ge0$ for every $U\subseteq\Omega$, and $\sum_{U\subseteq\Omega}\pi_x(U)=1$.
\item For every $i\in\Omega$, the sets containing $i$ have total weight $x_i$, and the sets avoiding $i$ have total weight $1-x_i$.
\item $\pi_x(U)>0$ if and only if $J_1\subseteq U\subseteq J_1\cup Q$. Consequently, $\pi_x\bigl(\binom\Omega j\bigr)>0$ if and only if $|J_1|\le j\le|J_1|+|Q|$.
\item If $\Omega=\Omega_1\sqcup\Omega_2$ and $x^{(1)},x^{(2)}$ denote the restrictions of $x$ to $\Omega_1,\Omega_2$, then $\pi_x(U)=\pi_{x^{(1)}}(U\cap\Omega_1)\,\pi_{x^{(2)}}(U\cap\Omega_2)$ for every $U\subseteq\Omega$.
\item $\pi_x(U)$ is a polynomial in the coordinates of $x$. In particular, it depends continuously on $x$.
\end{enumerate}
\end{lemma}
\begin{proof}
(a) Every factor in $\pi_x(U)$ lies in $[0,1]$. Expand the product
$\prod_{i\in\Omega}\bigl(x_i+(1-x_i)\bigr)=1$
by the distributive law, choosing one of the two summands from each factor. Each choice corresponds to exactly one set $U$, namely the set of indices at which $x_i$ was chosen, and the resulting term is $\pi_x(U)$. Hence the weights sum to $1$.

(b) Expand $x_i\prod_{\ell\ne i}\bigl(x_\ell+(1-x_\ell)\bigr)=x_i$ in the same way. The terms are exactly the weights $\pi_x(U)$ of the sets $U$ containing $i$. The statement about sets avoiding $i$ follows by subtracting from (a).

(c) A product of numbers in $[0,1]$ is positive exactly when every factor is positive. For $i\in U$ the factor $x_i$ is positive exactly when $i\notin J_0$; for $i\notin U$ the factor $1-x_i$ is positive exactly when $i\notin J_1$. This proves the first assertion. A set of size $j$ with positive weight is obtained by adding $j-|J_1|$ elements of $Q$ to $J_1$, which is possible exactly when $0\le j-|J_1|\le|Q|$.

(d) Regroup the factors of $\pi_x(U)$ according to whether their index lies in $\Omega_1$ or in $\Omega_2$.

(e) This is immediate from the definition.
\end{proof}

\subsection{Shares and averages on a layer}
\begin{definition}[Shares and averages]\label{def:shares}
Let $j$ be an integer with $\pi_x\bigl(\binom\Omega j\bigr)>0$. For $U\in\binom\Omega j$ and $\calF\subseteq\binom\Omega j$, define
$$
\nu_{x,j}(U)=\frac{\pi_x(U)}{\pi_x\bigl(\binom\Omega j\bigr)} \ \ \text{and } \
\nu_{x,j}(\calF)=\sum_{U\in\calF}\nu_{x,j}(U).
$$
We call $\nu_{x,j}(\calF)$ the \emph{share} of $\calF$ in the $j$-th layer. For a function $h:\binom\Omega j\to\R$, the \emph{weighted average} of $h$ is
$$
\avg{h}_{x,j}=\sum_{U\in\binom\Omega j}\nu_{x,j}(U)\,h(U).
$$
\end{definition}

Shares are nonnegative, additive over disjoint families, and the whole layer has share $1$. If every $x_i$ equals $1/2$, all sets of one layer have the same weight, so $\nu_{x,j}(\calF)=|\calF|/\binom{|\Omega|}j$ is the plain proportion and $\avg h_{x,j}$ is the plain average.

In the proof of Theorem~\ref{thm:criterion}, ``$\calF_c$ and $\calH_c$ are large'' will mean that both have share at least $1/n$ in the $k$-th layer of $P$, for one suitable vector $x$.

\subsection{The probabilistic reading}\label{subsec:dictionary}
The weights have a probabilistic meaning, which explains the terminology of the literature but is not needed for the proofs. Let $Z_x$ be the random subset of $\Omega$ obtained by including each $i$ independently with probability $x_i$.
\begin{center}
\begin{tabular}{@{}ll@{}}
\toprule
Weighted counting (this paper) & Probability\\
\midrule
$\pi_x(U)$ & $\Pr(Z_x=U)$\\
$\nu_{x,j}(\calF)$ & $\Pr(Z_x\in\calF\mid |Z_x|=j)$\\
$\avg h_{x,j}$ & expectation of $h(Z_x)$ given $|Z_x|=j$\\
nested transport (Definition~\ref{def:transport}) & coupling $(R,U)$ with $R\subseteq U$ almost surely\\
prefix averages $M_i$ (Lemma~\ref{lem:concentration}) & Doob martingale of $h(Z_x)$\\
the sums $\Phi_i(\lambda)$ (Lemma~\ref{lem:concentration}) & exponential moments of the martingale\\
Lemma~\ref{lem:two-point} & Hoeffding's lemma for a two-valued variable\\
\bottomrule
\end{tabular}
\end{center}
In this language, the coloring is used as follows: the color classes $\calF_i,\calG_i$ are fixed families determined by $\varphi$, and the random set $Z_x$ is drawn from a product distribution whose parameters $x=x^*$ are themselves determined by $\varphi$ (Section~\ref{sec:balance}). Nothing about $\varphi$ is random.

\subsection{Why the weights must depend on the coloring}
The uniform weights $x=(\tfrac12,\ldots,\tfrac12)$ cannot be used in general. The following example concerns arbitrary rules with $f(S)\in S$ and $g(R)\notin R$, not necessarily arising from a proper coloring.

\begin{remark}\label{rem:uniform}
Fix $X\in\binom Pk$, put $R_0=P\setminus X$, fix $z\in X$, and put $S_0=R_0\cup\{z\}$, so that $|R_0|=k-1$, $|S_0|=k$, and $X\cap S_0=\{z\}$.
\begin{itemize}
\item Every $k$-set $S$ meets $X$, because $|P\setminus X|=k-1$. Let $f(S)$ be an element of $S\cap X$ different from $z$ whenever such an element exists. This fails only if $S\cap X=\{z\}$, that is, only if $S=S_0$; put $f(S_0)=z$.
\item Every $(k-1)$-set $R$ misses an element of $S_0$, since $|S_0|=k$. Let $g(R)$ be an element of $S_0\setminus R$ different from $z$ whenever such an element exists. This fails only if $S_0\setminus R=\{z\}$, that is, only if $R=R_0$, let $g(R_0)=z$.
\end{itemize}
The rule $f$ uses colors in $X$ and $g$ uses colors in $S_0$, so $z$ is the only color used by both, and it is used only on $S_0$ and on $R_0$. With uniform weights, both classes of $z$ have share $1/\binom{2k-1}k$, which is exponentially small. The separation estimate of Section~\ref{sec:separation} applied with this share yields only $r\le\sqrt{2n\ln\binom{2k-1}k}$. The two middle coefficients $\binom{2k-1}{k-1}=\binom{2k-1}k$ are the largest of the $2k$ coefficients $\binom{2k-1}t$, which sum to $2^{2k-1}$; hence $\binom{2k-1}k\ge2^{2k-1}/(2k)$, and the bound is at least $\sqrt{2n(n\ln2-\ln(n+1))}$. This exceeds $n\ge k$ once $(2\ln2-1)n>2\ln(n+1)$, which holds for all large $k$. Since no two $k$-sets are at Johnson distance more than $k$, the bound then carries no information.

By contrast, the weights with $x_i=1$ on $R_0$, $x_i=0$ on $X\setminus\{z\}$, and $x_z=\tfrac12$ give positive weight only to $S_0$ in the $k$-th layer and only to $R_0$ in the $(k-1)$-th layer, so color $z$ has share $1$ on both sides. Section~\ref{sec:balance} finds such adapted weights for every coloring.
\end{remark}

\section{Choosing the weights from the coloring}\label{sec:balance}
We return to the standing assumption. This section is where the coloring determines the weights.

\subsection{Statement}
\begin{lemma}[Balanced weights]\label{lem:balance}
There exist $x^*\in[0,1]^n$ and a color $c\in P$ such that both $\pi_{x^*}\bigl(\binom Pk\bigr)$ and $\pi_{x^*}\bigl(\binom P{k-1}\bigr)$ are positive and
\begin{equation}\label{eq:balance}
\nu_{x^*,k}(\calF_c)=\nu_{x^*,k-1}(\calG_c)\ge\frac1n.
\end{equation}
\end{lemma}

The proof uses only that $\calF_1,\ldots,\calF_n$ and $\calG_1,\ldots,\calG_n$ partition the two layers and satisfy~\eqref{eq:class-containment}. It occupies the rest of this section.

\subsection{The set of admissible weights}\label{subsec:admissible}
Define
\begin{equation}\label{eq:D}
D=\left\{x\in[0,1]^n:\sum_{i=1}^n x_i=k-\tfrac12\right\}.
\end{equation}
Since $n/2=k-1/2$, the vector $(1/2,\ldots,1/2)$ belongs to $D$, so $D$ is nonempty. It is the intersection of a closed cube with a closed affine hyperplane. Consequently, it is closed and bounded in $\R^n$, and hence compact. The cube and the hyperplane are both convex, so their intersection $D$ is convex as well.

For $x\in D$, define
$a(x)=\pi_x\bigl(\tbinom Pk\bigr)$ and $b(x)=\pi_x\bigl(\tbinom P{k-1}\bigr)
$
for the total weights of the two layers. We show that both are positive for every $x\in D$, including points with coordinates equal to $0$ or $1$. Let $u$ be the number of coordinates with $x_i=1$, and let $v$ be the number of coordinates with $x_i>0$. Since
\[
u\le\sum_{i=1}^n x_i=k-\tfrac12\le v
\]
and $u,v$ are integers, we have $u\le k-1$ and $v\ge k$. By Lemma~\ref{lem:weights}(c), the $j$-th layer has positive weight exactly when $u\le j\le v$. In particular,
\begin{equation}\label{eq:ab-positive}
a(x)>0 \ \ \text{and } \ b(x)>0\qquad(x\in D).
\end{equation}
The shares $\nu_{x,k}$ and $\nu_{x,k-1}$ are therefore defined throughout $D$. This is the reason for the constraint $\sum_ix_i=k-\tfrac12$: it places the ``typical size'' of a set halfway between the two layers that matter.

\subsection{A continuous self-map built from the coloring}\label{subsec:selfmap}
 Let
$A_i(x)=\pi_x(\calF_i)$ and $B_i(x)=\pi_x(\calG_i)$ for each $i\in P$ and $x\in D$.
These are the unnormalized weights of the color classes of $i$; this is the point at which the coloring enters. They are finite sums of polynomials in $x$, so every $A_i$ and $B_i$ is continuous on $D$ by Lemma~\ref{lem:weights}(e). Since the classes partition the two layers,
\begin{equation}\label{eq:sumAB}
\sum_{i=1}^n A_i(x)=a(x),\qquad
\sum_{i=1}^n B_i(x)=b(x).
\end{equation}
By~\eqref{eq:class-containment} and Lemma~\ref{lem:weights}(b), we have
\begin{equation}\label{eq:AB-bound}
0\le A_i(x)\le x_i \ \ \text{and } \ 0\le B_i(x)\le1-x_i.
\end{equation}
Also $0\le a(x),b(x)\le1$ by Lemma~\ref{lem:weights}(a).

Define a map $T:D\to\R^n$ by
\begin{equation}\label{eq:T}
T_i(x)=x_i+a(x)B_i(x)-b(x)A_i(x)
\qquad(1\le i\le n).
\end{equation}
Each coordinate is continuous, and therefore $T$ is continuous. We verify that $T(x)\in D$ whenever $x\in D$. Using~\eqref{eq:AB-bound} and $0\le a(x),b(x)\le1$,
\begin{align*}
T_i(x)&\ge x_i-b(x)A_i(x)\ge x_i-A_i(x)\ge0,\\
T_i(x)&\le x_i+a(x)B_i(x)\le x_i+B_i(x)\le1.
\end{align*}
Thus $T(x)\in[0,1]^n$. Summing~\eqref{eq:T} over $i$ and using~\eqref{eq:sumAB},
\[
\sum_{i=1}^n T_i(x)
=\sum_{i=1}^n x_i+a(x)b(x)-b(x)a(x)
=k-\tfrac12.
\]
This proves $T(D)\subseteq D$.

The map has the following heuristic meaning, which is not used in the proof. Dividing the correction term by $a(x)b(x)$, it equals the share of color $i$ under $g$ minus its share under $f$. If $g$ uses color $i$ more heavily than $f$ does, then $T$ raises $x_i$. This makes the sets containing $i$ heavier, which tends to increase the share of $\calF_i$ and to decrease the share of $\calG_i$.

\subsection{The fixed point}
We use Brouwer's fixed-point theorem~\cite{Brouwer1911}: every continuous map from a closed simplex in $\R^q$ to itself has a fixed point. Here a closed simplex is the convex hull of finitely many affinely independent points. In dimension one the theorem is the intermediate value theorem: if $F:[0,1]\to[0,1]$ is continuous, then $F(x)-x$ is $\ge0$ at $0$ and $\le0$ at $1$, so it vanishes somewhere.

The theorem does not apply to $T:D\to D$ directly, because $D$ is not a simplex: it is a polytope of dimension $n-1$ lying in a hyperplane of $\R^n$. We therefore enlarge $D$ to a simplex $\Sigma$, and pull $\Sigma$ back onto $D$ by a continuous map that fixes every point of $D$.

\emph{The simplex.} Let
\[
\Sigma=\Bigl\{z\in\R^n:\ z_i\ge0\ (1\le i\le n),\ \sum_{i=1}^nz_i\le n\Bigr\}.
\]
This is the convex hull of the $n+1$ affinely independent points $0,ne_1,\ldots,ne_n$, where $e_1,\ldots,e_n$ are the unit vectors. Indeed, a point $z$ with $z\ge0$ and $\sum_iz_i\le n$ equals $\sum_i\frac{z_i}n(ne_i)+\bigl(1-\sum_i\frac{z_i}n\bigr)\cdot0$, a convex combination of these points; conversely, every convex combination of them is nonnegative with coordinate sum at most $n$. Every $x\in D$ satisfies $0\le x_i\le1$, hence $x\ge0$ and $\sum_ix_i\le n$. Thus $D\subseteq\Sigma$.

\emph{The nearest-point map.} We write $\langle\cdot,\cdot\rangle$ and $\|\cdot\|$ for the Euclidean inner product and norm on $\R^n$.

\begin{lemma}[Nearest points of $D$]\label{lem:nearest}
For every $z\in\R^n$ there is exactly one point $p_D(z)\in D$ with $\|z-p_D(z)\|=\min_{y\in D}\|z-y\|$. Moreover:
\begin{enumerate}[label=\textup{(\alph*)}]
\item $\langle z-p_D(z),\,y-p_D(z)\rangle\le0$ for every $y\in D$;
\item $\|p_D(z)-p_D(z')\|\le\|z-z'\|$ for all $z,z'\in\R^n$; in particular, $p_D$ is continuous;
\item $p_D(z)=z$ for every $z\in D$.
\end{enumerate}
\end{lemma}
\begin{proof}
By Subsection~\ref{subsec:admissible}, $D$ is nonempty, compact, and convex.

\emph{Existence.} The function $y\mapsto\|z-y\|$ is continuous on the nonempty compact set $D$, so it attains its minimum at some point $p$.

\emph{Inequality (a) for every minimizer $p$.} Let $y\in D$ and $0<t\le1$. By convexity, $p+t(y-p)\in D$, so by minimality
$\|z-p-t(y-p)\|^2\ge\|z-p\|^2$.
Expanding the left side gives $-2t\langle z-p,y-p\rangle+t^2\|y-p\|^2\ge0$. Dividing by $t$ and letting $t\to0$ yields $\langle z-p,y-p\rangle\le0$.

\emph{Uniqueness.} Let $p$ and $p'$ be minimizers. Applying (a) at $p$ with $y=p'$ and at $p'$ with $y=p$ gives
$
\langle z-p,\,p'-p\rangle\le0$ and $\langle z-p',\,p-p'\rangle\le0$.
Adding the two inequalities gives $\|p'-p\|^2\le0$, so $p'=p$. Hence $p_D(z)$ is well defined, and (a) holds for it.

\emph{(b).} Let $p=p_D(z)$ and $p'=p_D(z')$. By (a), $\langle z-p,\,p'-p\rangle\le0$ and $\langle z'-p',\,p-p'\rangle\le0$. Adding these and rearranging,
$\|p-p'\|^2\le\langle z-z',\,p-p'\rangle\le\|z-z'\|\,\|p-p'\|$,
by the Cauchy--Schwarz inequality. If $p\ne p'$, dividing by $\|p-p'\|$ gives (b); if $p=p'$, (b) holds trivially.

\emph{(c).} If $z\in D$, then $y=z$ gives distance $0$, which is the least possible value, so $p_D(z)=z$ by uniqueness.
\end{proof}

\emph{Applying Brouwer's theorem.} Define $h:\Sigma\to\R^n$ by $h(z)=T(p_D(z))$. By Lemma~\ref{lem:nearest}(b), $p_D$ is continuous, and $T$ is continuous on $D$ by Subsection~\ref{subsec:selfmap}; hence $h$ is continuous. Moreover $p_D(z)\in D$, $T(D)\subseteq D$ by Subsection~\ref{subsec:selfmap}, and $D\subseteq\Sigma$; hence $h(\Sigma)\subseteq D\subseteq\Sigma$. Brouwer's theorem therefore gives a point $z\in\Sigma$ with $T(p_D(z))=z$.
The left side lies in $D$, hence $z\in D$. By Lemma~\ref{lem:nearest}(c), $p_D(z)=z$, and substituting into $T(p_D(z))=z$ gives $T(z)=z$. Thus there exists $x^*\in D$ with
\begin{equation}\label{eq:fixed}
T(x^*)=x^*.
\end{equation}

\begin{proof}[Proof of Lemma~\ref{lem:balance}]
Let $x^*$ be as in~\eqref{eq:fixed}. Both layers have positive $x^*$-weight by~\eqref{eq:ab-positive}.
Substituting $T(x^*)=x^*$ into~\eqref{eq:T} yields $a(x^*)B_i(x^*)=b(x^*)A_i(x^*)$ for each $i\in P$. Dividing by $a(x^*)b(x^*)>0$ gives
\begin{equation}\label{eq:equal-output}
\nu_{x^*,k}(\calF_i)=\frac{A_i(x^*)}{a(x^*)}=\frac{B_i(x^*)}{b(x^*)}=\nu_{x^*,k-1}(\calG_i)=:q_i.
\end{equation}
Thus, under the weights $x^*$, every color has the same share as an $f$-color among the $k$-sets as it has as a $g$-color among the $(k-1)$-sets. By~\eqref{eq:sumAB}, the numbers $q_i$ are nonnegative and sum to $1$. If every $q_i$ were smaller than $1/n$, their sum would be smaller than $1$. Consequently, there is a color $c\in P$ such that $q_c\ge1/n$.
\end{proof}

From now on, $x^*$ and $c$ are fixed as in Lemma~\ref{lem:balance}.

\section{Lifting the second class to the upper layer}\label{sec:lift}
Lemma~\ref{lem:balance} makes $\calF_c$ large among $k$-sets and $\calG_c$ large among $(k-1)$-sets. Lemma~\ref{lem:far} compares $\calF_c$ with $\calH_c=\nabla\calG_c$. This section shows that $\calH_c$ is at least as large as $\calG_c$:
\begin{equation}\label{eq:lift}
\nu_{x^*,k}(\calH_c)\ge\nu_{x^*,k-1}(\calG_c)\ge\frac1n.
\end{equation}
This is not automatic: a single $k$-set contains $k$ different $(k-1)$-sets, and with nonuniform weights a naive count can lose weight. The tool is a way of shipping the weight of each $(k-1)$-set to its supersets.

\subsection{Nested transports and the upward transfer}
\begin{definition}[Nested transport]\label{def:transport}
Let $\Omega$ be a finite set, and $j$ an integer. Let $\alpha$ and $\beta$ be nonnegative functions of $\binom\Omega{j-1}$ and $\binom\Omega j$, each with a total of $1$. A \emph{nested transport} from $\alpha$ to $\beta$ is an assignment of a number $t(R,U)\ge0$ to every pair $R\in\binom\Omega{j-1}$, $U\in\binom\Omega j$ with $R\subseteq U$, such that
\begin{equation}\label{eq:transport}
\sum_{U\supseteq R}t(R,U)=\alpha(R)\quad\Bigl(R\in\tbinom\Omega{j-1}\Bigr)
\ \ \text{and } \
\sum_{R\subseteq U}t(R,U)=\beta(U)\quad\Bigl(U\in\tbinom\Omega{j}\Bigr).
\end{equation}
Here the sums run over $U\in\binom\Omega j$ and $R\in\binom\Omega{j-1}$, respectively.
\end{definition}

In words, each set $R$ of the lower layer ships its amount $\alpha(R)$ to its supersets, and each set $U$ of the upper layer receives exactly $\beta(U)$.

\begin{example}\label{ex:uniform-transport}
Let $|\Omega|=2k-1$, let $\alpha$ and $\beta$ be the uniform functions on the layers $k-1$ and $k$, and put $M=\binom{2k-1}{k}=\binom{2k-1}{k-1}$. Every $(k-1)$-set has exactly $k$ supersets in the $k$-th layer, and every $k$-set has exactly $k$ subsets in the $(k-1)$-th layer. Hence $t(R,U)=1/(kM)$ for all $R\subseteq U$ is a nested transport. In general the weights are not uniform, and the construction below is needed.
\end{example}

\begin{lemma}[Upward transfer]\label{lem:shadow}
If there is a nested transport from $\alpha$ to $\beta$, then $\beta(\nabla\calG)\ge\alpha(\calG)$ for every $\calG\subseteq\binom\Omega{j-1}$.
\end{lemma}
\begin{proof}
By the first identity in~\eqref{eq:transport},
$$
\alpha(\calG)=\sum_{R\in\calG}\ \sum_{U\supseteq R}t(R,U).
$$
Every pair $(R,U)$ in this double sum has $U\in\nabla\calG$. Adding the nonnegative terms $t(R,U)$ with $U\in\nabla\calG$ and $R\notin\calG$, and then applying the second identity in~\eqref{eq:transport}, we obtain
$$
\alpha(\calG)\le\sum_{U\in\nabla\calG}\ \sum_{R\subseteq U}t(R,U)
=\beta(\nabla\calG).\qedhere
$$
\end{proof}

\subsection{Existence of nested transports, and the lift}
\begin{proposition}[Nested transport for product weights]\label{prop:transport}
Let $\Omega$ be a finite set, $x\in[0,1]^\Omega$, and $j$ an integer. If $\pi_x\bigl(\binom\Omega{j-1}\bigr)>0$ and $\pi_x\bigl(\binom\Omega{j}\bigr)>0$, then there is a nested transport from $\nu_{x,j-1}$ to $\nu_{x,j}$. Consequently, $\nu_{x,j}(\nabla\calG)\ge\nu_{x,j-1}(\calG)$ with $\Bigl(\calG\subseteq\tbinom\Omega{j-1}\Bigr)$.

\end{proposition}

\begin{corollary}\label{cor:lift}
Inequality~\eqref{eq:lift} holds.
\end{corollary}
\begin{proof}
Apply Proposition~\ref{prop:transport} with $\Omega=P$, $x=x^*$, $j=k$ and $\calG=\calG_c$, both layers have positive weight according to Lemma~\ref{lem:balance}. Then use~\eqref{eq:balance}.
\end{proof}

The proof of Proposition~\ref{prop:transport} occupies the rest of this section. Related statements belong to the theory of strong Rayleigh measures; see Pemantle and Peres~\cite{PP2014}. We give direct proofs of precisely the statements required here and do not claim novelty for them.

\subsection{Positive weights}
Let $w_1,\ldots,w_N>0$. Let $w(U)=\prod_{i\in U}w_i$ for $U\subseteq[N]$, and for $0\le j\le N$ define the total layer and the normalized weight
\begin{equation}\label{eq:weighted}
W_j=W_j(w_1,\ldots,w_N)=\sum_{U\in\binom{[N]}j}w(U)>0,
\qquad
\omega_j(U)=\frac{w(U)}{W_j}\quad\Bigl(U\in\tbinom{[N]}j\Bigr).
\end{equation}
We set $W_j=0$ for $j<0$ and $j>N$. Expanding the product by choosing either $1$ or $w_i z$ from each factor gives
\begin{equation}\label{eq:generating}
\prod_{i=1}^N(1+w_i z)=\sum_{j=0}^N W_jz^j.
\end{equation}

\begin{lemma}\label{lem:logconcave}
The coefficients in~\eqref{eq:generating} form a log-concave sequence: for every integer $j$,
$W_j^2\ge W_{j-1}W_{j+1}$.
\end{lemma}
\begin{proof}
We first establish a preservation property. Let $(a_j)_{j\in\mathbb Z}$ be a nonzero, finitely supported sequence of nonnegative numbers whose positive entries have consecutive indices. Suppose that
$a_j^2\ge a_{j-1}a_{j+1}$ with $j\in\mathbb Z$.
Given $w>0$, let $b_j=a_j+wa_{j-1}$. Expanding both products and collecting the coefficients of $1,w,w^2$, we obtain
\begin{align}\label{eq:lc-expansion}
b_j^2-b_{j-1}b_{j+1}=a_j^2-a_{j-1}a_{j+1}+w(a_ja_{j-1}-a_{j-2}a_{j+1})+w^2(a_{j-1}^2-a_{j-2}a_j).
\end{align}
The first and third terms are nonnegative by the assumed log-concavity. To treat the middle term, first suppose that $a_{j-2}a_{j+1}=0$. Then its nonnegativity follows from $a_ja_{j-1}\ge0$. If instead $a_{j-2}a_{j+1}>0$, consecutiveness of the positive support implies that all four entries $a_{j-2},a_{j-1},a_j,a_{j+1}$ are positive. Dividing the log-concavity inequalities at indices $j-1$ and $j$ by the relevant positive products gives
\[
\frac{a_{j-1}}{a_{j-2}}
\ge\frac{a_j}{a_{j-1}}
\ge\frac{a_{j+1}}{a_j}.
\]
Comparing the first and last ratios and multiplying by $a_{j-2}a_j>0$ yields
$a_{j-1}a_j\ge a_{j-2}a_{j+1}$.
Thus every term on the right-hand side of~\eqref{eq:lc-expansion} is nonnegative, and $(b_j)$ is log-concave. Moreover, if the positive support of $(a_j)$ is $\{\ell,\ldots,u\}$, the positive support of $(b_j)$ is exactly $\{\ell,\ldots,u+1\}$: the two shifted supports overlap or are consecutive, and $w>0$. Hence the support condition is preserved as well.

Start with $a_0=1$ and $a_j=0$ for $j\ne0$. This sequence is log-concave and has consecutive positive support. Successively multiply its generating polynomial by $1+w_1z,\ldots,1+w_Nz$. Multiplication by $1+wz$ replaces the coefficient sequence $(a_j)$ by $(a_j+wa_{j-1})$, so the preceding argument applies at every step. By~\eqref{eq:generating}, the final coefficient sequence is $(W_j)$, proving the lemma.
\end{proof}

\begin{lemma}[Nested transport for positive weights]\label{lem:coupling}
For positive weights $w_1,\ldots,w_N$ and every $1\le j\le N$, there is a nested transport from $\omega_{j-1}$ to $\omega_j$.
\end{lemma}
\begin{proof}
We argue by induction on $N$. For $N=1$, necessarily $j=1$; the layers are $\{\varnothing\}$ and $\{\{1\}\}$, both normalized weights equal $1$, and $t(\varnothing,\{1\})=1$ is a nested transport.

Suppose now that $N\ge2$ and that the assertion holds for $N-1$ elements. Write $W'_t$ and $\omega'_t$ for the quantities~\eqref{eq:weighted} computed from $w_1,\ldots,w_{N-1}$; thus $W'_t>0$ for $0\le t\le N-1$ and $W'_t=0$ otherwise. Splitting the sets of size $t$ according to whether they contain $N$ gives
\begin{equation}\label{eq:split}
W_t=W'_t+w_NW'_{t-1}.
\end{equation}
For $0\le t\le N$, let
\begin{equation}\label{eq:pj}
p_t=\frac{w_NW'_{t-1}}{W_t}
\end{equation}
be the share of the $t$-th layer of $[N]$ occupied by sets containing $N$. Then $p_0=0$ and $p_N=1$, and $1-p_t=W'_t/W_t$ by~\eqref{eq:split}. Directly from the definitions,
\begin{align}
\omega_t(V)&=(1-p_t)\,\omega'_t(V)
&&\bigl(V\subseteq[N-1],\ |V|=t,\ 0\le t\le N-1\bigr),\label{eq:I1}\\
\omega_t(V\cup\{N\})&=p_t\,\omega'_{t-1}(V)
&&\bigl(V\subseteq[N-1],\ |V|=t-1,\ 1\le t\le N\bigr).\label{eq:I2}
\end{align}
Indeed, the right-hand side of~\eqref{eq:I1} is $(W'_t/W_t)(w(V)/W'_t)=w(V)/W_t$, and that of~\eqref{eq:I2} is $(w_NW'_{t-1}/W_t)(w(V)/W'_{t-1})=w(V\cup\{N\})/W_t$; all denominators are positive in the stated ranges.

We next show $p_{j-1}\le p_j$ for $1\le j\le N$. For $j=1$ this reads $0\le p_1$, and for $j=N$ it reads $p_{N-1}\le1$. For $2\le j\le N-1$, subtraction over a common denominator yields
$$
p_j-p_{j-1}
=\frac{w_N\big((W'_{j-1})^2-W'_{j-2}W'_j\big)}
{W_jW_{j-1}}
\ge0,
$$
because the denominator is positive and the numerator is nonnegative by Lemma~\ref{lem:logconcave} applied to $w_1,\ldots,w_{N-1}$.
Fix $1\le j\le N$. Every pair $R\subseteq U$ with $|R|=j-1$ and $|U|=j$ is of exactly one of the following three kinds, and we define $t(R,U)$ accordingly.
\begin{enumerate}[label=\textup{(K\arabic*)}]
\item $N\notin U$. Such pairs exist only if $j\le N-1$; then the induction hypothesis gives a nested transport $t'$ from $\omega'_{j-1}$ to $\omega'_j$ on $[N-1]$, and we put $t(R,U)=(1-p_j)\,t'(R,U)$.
\item $N\notin R$ and $N\in U$. Then $U=R\cup\{N\}$, since $U\setminus R$ has one element. Put $t(R,U)=(p_j-p_{j-1})\,\omega'_{j-1}(R)$.
\item $N\in R$. Write $R=R'\cup\{N\}$ and $U=U'\cup\{N\}$, so $R'\subseteq U'\subseteq[N-1]$ and $|R'|=j-2$. Such pairs exist only if $j\ge2$; then $1\le j-1\le N-1$, the induction hypothesis gives a nested transport $t''$ from $\omega'_{j-2}$ to $\omega'_{j-1}$ on $[N-1]$, and we put $t(R,U)=p_{j-1}\,t''(R',U')$.
\end{enumerate}
All these numbers are nonnegative because $0\le p_{j-1}\le p_j\le1$.

\emph{Amount shipped by $R$ with $N\notin R$.} The supersets of $R$ are the sets $U\subseteq[N-1]$ of kind (K1), together with $R\cup\{N\}$ of kind (K2). If $j\le N-1$, the first group ships $(1-p_j)\omega'_{j-1}(R)$ in total; if $j=N$, there is no such $U$, and $1-p_N=0$, so the same formula holds. Adding the amount of kind (K2) gives
$
(1-p_j)\omega'_{j-1}(R)+(p_j-p_{j-1})\omega'_{j-1}(R)
=(1-p_{j-1})\omega'_{j-1}(R)=\omega_{j-1}(R)$
by~\eqref{eq:I1} with $t=j-1$.

\emph{Amount shipped by $R=R'\cup\{N\}$.} This requires $j\ge2$. All supersets are of kind (K3), and the total is $p_{j-1}\omega'_{j-2}(R')=\omega_{j-1}(R)$ by~\eqref{eq:I2} with $t=j-1$.

\emph{Amount received by $U$ with $N\notin U$.} This requires $j\le N-1$. All subsets are of kind (K1), and the total is $(1-p_j)\omega'_j(U)=\omega_j(U)$ by~\eqref{eq:I1} with $t=j$.

\emph{Amount received by $U=U'\cup\{N\}$.} The subsets of $U$ in the $(j-1)$-th layer are $U'$, of kind (K2), and the sets $R'\cup\{N\}$ with $R'\subseteq U'$, of kind (K3). The total is
$
(p_j-p_{j-1})\omega'_{j-1}(U')+p_{j-1}\omega'_{j-1}(U')=p_j\omega'_{j-1}(U')=\omega_j(U)$
by~\eqref{eq:I2} with $t=j$. When $j=1$ there are no sets of kind (K3), and the formula remains correct because $p_0=0$.

Thus $t$ satisfies~\eqref{eq:transport}, completing the induction.
\end{proof}

\subsection{Forced and forbidden elements}
\begin{lemma}[Reduction to positive weights]\label{lem:reduction}
Let $\Omega$ be a finite set and $x\in[0,1]^\Omega$. Let $J_0,J_1,Q$ be as in Lemma~\ref{lem:weights}, put $u=|J_1|$, and give each $i\in Q$ the positive weight $w_i=x_i/(1-x_i)$. If $\pi_x\bigl(\binom\Omega j\bigr)>0$, then for every $U\in\binom\Omega j$,
\[
\nu_{x,j}(U)=
\begin{cases}
\omega_{j-u}(U\setminus J_1) & \text{if }J_1\subseteq U\subseteq J_1\cup Q,\\
0&\text{otherwise},
\end{cases}
\]
where $\omega_{j-u}$ is computed on the ground set $Q$ from the weights $(w_i)_{i\in Q}$.
\end{lemma}
\begin{proof}
If $U$ is not of the form $J_1\cup V$ with $V\subseteq Q$, then $\pi_x(U)=0$ by Lemma~\ref{lem:weights}(c). If $U=J_1\cup V$ with $V\subseteq Q$, then the factors for $i\in J_0\cup J_1$ all equal $1$, and
\[
\pi_x(U)=\prod_{i\in V}x_i\prod_{i\in Q\setminus V}(1-x_i)
=\gamma\prod_{i\in V}w_i,
\qquad \gamma=\prod_{i\in Q}(1-x_i)>0.
\]
Summing over all $V\in\binom Q{j-u}$ gives $\pi_x\bigl(\binom\Omega j\bigr)=\gamma\,W_{j-u}$, and dividing proves the formula.
\end{proof}

\begin{proof}[Proof of Proposition~\ref{prop:transport}]
With the notation of Lemma~\ref{lem:reduction}, Lemma~\ref{lem:weights}(c) applied to both layers gives $u\le j-1$ and $j\le u+|Q|$, that is, $1\le j-u\le|Q|$. In particular $Q\ne\varnothing$. Label the elements of $Q$ by $1,\ldots,|Q|$ and apply Lemma~\ref{lem:coupling} to obtain a nested transport $t_Q$ from $\omega_{j-1-u}$ to $\omega_{j-u}$ on $Q$. Let
$t(J_1\cup V,\,J_1\cup V')=t_Q(V,V')$ when $V\subseteq V'\subseteq Q$,
and $t(R,U)=0$ for all other pairs $R\subseteq U$. Each pair receiving a possibly positive value satisfies $J_1\cup V\subseteq J_1\cup V'$. By Lemma~\ref{lem:reduction}, a set $R$ of the form $J_1\cup V$ ships $\sum_{V'}t_Q(V,V')=\omega_{j-1-u}(V)=\nu_{x,j-1}(R)$, and every other set of the lower layer ships $0=\nu_{x,j-1}(R)$. The amounts received are checked in the same way. The final inequality is Lemma~\ref{lem:shadow}.
\end{proof}

\section{Two large families on one layer are close}\label{sec:separation}
By Lemma~\ref{lem:balance} and Corollary~\ref{cor:lift}, the families $\calF_c$ and $\calH_c$ both have share at least $1/n$ in the $k$-th layer under the weights $x^*$. This section shows that two such families cannot be far apart. Its results are general facts about weights and do not use $\varphi$.

\subsection{Statement}
\begin{lemma}[Separation of two families]\label{lem:separation}
Let $\Omega=[N]$ with $N\ge1$, let $x\in[0,1]^N$, let $\pi_x\bigl(\binom{[N]}j\bigr)>0$, and write $\nu=\nu_{x,j}$. Suppose that
$\calF,\calH\subseteq\tbinom{[N]}j$,
$\nu(\calF),\nu(\calH)\ge\rho$ and $0<\rho\le1$.
If $r\ge0$ satisfies $d_J(S,U)\ge r$ for every $S\in\calF$ and $U\in\calH$, then
\begin{equation}\label{eq:separation}
r\le\sqrt{2N\ln(1/\rho)}.
\end{equation}
\end{lemma}

The idea is an isoperimetric one. Let $h(U)$ be the distance from $U$ to $\calF$. A single exchange changes $h$ by at most $1$, and such a function is nearly constant on most of the layer (Lemma~\ref{lem:concentration}). Since $h$ vanishes on the large family $\calF$, its average is small; since $h\ge r$ on the large family $\calH$, its average is large. Comparing the two gives~\eqref{eq:separation}.

\subsection{Slowly changing functions}
Sets at Johnson distance one differ by a single exchange: one element is removed and another is inserted. Conversely, if $d_J(U,W)=q$, enumerate the differences as
$U\setminus W=\{u_1,\ldots,u_q\}$ and $W\setminus U=\{w_1,\ldots,w_q\}$.
Starting at $U$, replace $u_s$ by $w_s$ for $s=1,\ldots,q$. Each intermediate set has size $j$, each step is one exchange, and the final set is $W$. A real-valued function $h$ on a layer is \emph{$1$-Lipschitz} if
$|h(U)-h(W)|\le d_J(U,W)$ when $\left(U,W\in\tbinom{\Omega}j\right)$.
It suffices to check this inequality when $d_J(U,W)=1$: summing the resulting bounds along the exchange sequence above proves the general case.

\begin{lemma}[Concentration on a layer]\label{lem:concentration}
Let $\Omega=[N]$ with $N\ge1$, let $x\in[0,1]^N$, and let $j$ be an integer with $\pi_x\bigl(\binom{[N]}j\bigr)>0$. Write $\nu=\nu_{x,j}$. If $h:\binom{[N]}j\to\R$ is $1$-Lipschitz and $\bar h=\avg{h}_{x,j}$, then, for every $t\ge0$,
\begin{align}
\nu\bigl(\{U:h(U)\ge\bar h+t\}\bigr)&\le e^{-2t^2/N},\label{eq:upper}\\
\nu\bigl(\{U:h(U)\le\bar h-t\}\bigr)&\le e^{-2t^2/N}.\label{eq:lower}
\end{align}
\end{lemma}

\subsection{Proof of Lemma~\ref{lem:separation} from Lemma~\ref{lem:concentration}}
\begin{proof}[Proof of Lemma~\ref{lem:separation}]
Because $\nu(\calF)>0$, the family $\calF$ is nonempty. Define its distance function on the entire layer by
$h(U)=\min_{S\in\calF}d_J(U,S)$.
The minimum exists because $\calF$ is finite. To verify the Lipschitz property, fix $U,W$ and choose $S_W\in\calF$ with $h(W)=d_J(W,S_W)$. The triangle inequality gives
$
h(U)\le d_J(U,S_W)
\le d_J(U,W)+d_J(W,S_W)
=d_J(U,W)+h(W)$.
Interchanging $U,W$ yields $|h(U)-h(W)|\le d_J(U,W)$. Thus $h$ is $1$-Lipschitz. It is also nonnegative, and $h(U)=0$ if and only if $U\in\calF$.

Let $\eta=\avg{h}_{x,j}\ge0$ and
$\tau=\sqrt{\frac{N\ln(1/\rho)}2}$.
Since $h\ge0$, the families $\{U:h(U)=0\}$ and $\{U:h(U)\le\eta-\eta\}$ coincide. Apply~\eqref{eq:lower} with $t=\eta$ to obtain
$\rho\le\nu(\calF)=\nu\bigl(\{U:h(U)=0\}\bigr)
\le e^{-2\eta^2/N}$.
Taking natural logarithms gives $\ln\rho\le-2\eta^2/N$, or equivalently
$\eta^2\le\frac N2\ln(1/\rho)=\tau^2$.
As both $\eta$ and $\tau$ are nonnegative, we conclude that $\eta\le\tau$.

If $r\le\eta$, then $r\le\tau\le2\tau$, which already proves the required bound. If $r>\eta$, the separation hypothesis implies $h(U)\ge r$ for every $U\in\calH$: every distance in the minimum defining $h(U)$ is at least $r$. Applying~\eqref{eq:upper} with the positive threshold $t=r-\eta$ yields
$
\rho\le\nu(\calH)
\le\nu\bigl(\{U:h(U)\ge r\}\bigr)
\le e^{-2(r-\eta)^2/N}$.
Taking logarithms and using $r-\eta>0$ gives $r-\eta\le\tau$. Combining this with $\eta\le\tau$ proves
$r\le\eta+\tau\le2\tau
=\sqrt{2N\ln(1/\rho)}$.

This also covers $\rho=1$. Indeed, then $\tau=0$ and the first estimate forces $\eta=0$. If $r>0$, the second estimate would give $1\le e^{-2r^2/N}<1$, a contradiction. Hence $r=0$ in this boundary case, as asserted by~\eqref{eq:separation}.
\end{proof}

\subsection{A two-point inequality}
The only analytic input of the concentration lemma is the following one-variable inequality.

\begin{lemma}[A two-point inequality]\label{lem:two-point}
For every $p\in[0,1]$ and every real $s$,
$pe^{(1-p)s}+(1-p)\,e^{-ps}\le e^{s^2/8}$.
\end{lemma}
\begin{proof}
Put $D(s)=1-p+pe^{s}$. This is positive for every $s$: if $p<1$ then $1-p>0$, and if $p=1$ then $D(s)=e^s$. The left-hand side equals $e^{-ps}D(s)$, hence its logarithm is
$\gamma(s)=-ps+\ln D(s)$.
We have $\gamma(0)=\ln1=0$ and
$\gamma'(s)=-p+\frac{pe^{s}}{D(s)}$ with $\gamma'(0)=0$.
Differentiating once more and writing $a=pe^{s}\ge0$ and $b=1-p\ge0$, so that $a+b=D(s)>0$,
\[
\gamma''(s)=\frac{pe^{s}D(s)-(pe^{s})^2}{D(s)^2}=\frac{ab}{(a+b)^2}\le\frac14,
\]
where the last step is $(a+b)^2-4ab=(a-b)^2\ge0$. Let $H(s)=s^2/8-\gamma(s)$. Then $H(0)=H'(0)=0$ and $H''(s)=1/4-\gamma''(s)\ge0$. By the mean value theorem, $H'$ is nondecreasing, so $H'\ge0$ on $[0,\infty)$ and $H'\le0$ on $(-\infty,0]$. Applying the mean value theorem again, $H$ is nondecreasing on $[0,\infty)$ and nonincreasing on $(-\infty,0]$, so $H(s)\ge H(0)=0$ for every $s$. Exponentiating $\gamma(s)\le s^2/8$ proves the lemma.
\end{proof}

\subsection{Proof of the concentration lemma}
The idea is to decide, one element at a time, whether it belongs to the set. Each decision changes the relevant weighted average by at most $1$, and $N$ such changes rarely add up to much more than $\sqrt N$.

\begin{proof}[Proof of Lemma~\ref{lem:concentration}]
Let $\calS=\{U\in\binom{[N]}j:\pi_x(U)>0\}$, so that $\nu(\calS)=1$. For $0\le i\le N$ and $U\in\calS$, let
$C_i(U)=\{U'\in\calS:\ U'\cap[i]=U\cap[i]\}$
be the \emph{class} of $U$ at level $i$, consisting of the sets of positive weight that agree with $U$ on the first $i$ elements. Since $U\in C_i(U)$, we have $\nu(C_i(U))>0$, and we may define the \emph{prefix average}
\[
M_i(U)=\frac{1}{\nu(C_i(U))}\sum_{U'\in C_i(U)}\nu(U')\,h(U').
\]
The classes at level $i$ partition $\calS$, and $M_i$ is constant on each of them. At level $0$ there is a single class $\calS$, so $M_0(U)=\bar h$; at level $N$ every class is a single set, so $M_N(U)=h(U)$.

Fix $1\le i\le N$ and a class $C$ at level $i-1$. Its members share the trace $Y=U\cap[i-1]$; put $q=|Y|$ and let $m$ be the value of $M_{i-1}$ on $C$. The class splits into $C^1=\{U\in C:\ i\in U\}$ and $C^0=\{U\in C:\ i\notin U\}$,
each of which is empty or a class at level $i$. Let $u_1$ and $u_0$ be the values of $M_i$ on $C^1$ and $C^0$, when these are nonempty. Since the sum defining $M_{i-1}$ on $C$ splits along $C=C^1\sqcup C^0$,
\begin{equation}\label{eq:split-average}
\nu(C)\,m=\nu(C^1)\,u_1+\nu(C^0)\,u_0,
\end{equation}
where a term is omitted if its class is empty.

\emph{Step 1: one decision moves the average by at most one.} Suppose that $C^1$ and $C^0$ are both nonempty. We claim that $|u_1-u_0|\le1$.

Let $Q_i=\{i+1,\ldots,N\}$ and let $y$ be the restriction of $x$ to $Q_i$. Fix $\varepsilon\in\{0,1\}$ and put $Z=\{i\}$ if $\varepsilon=1$ and $Z=\varnothing$ if $\varepsilon=0$. Every set in the $j$-th layer that agrees with $Y$ on $[i-1]$ and satisfies the $\varepsilon$-condition at $i$ has the form $Y\cup Z\cup V$ with $V\in\binom{Q_i}{j-q-\varepsilon}$, and every such $V$ gives such a set. By Lemma~\ref{lem:weights}(d) with $\Omega=[i]\sqcup Q_i$, $\pi_x(Y\cup Z\cup V)=\kappa_\varepsilon\,\pi_y(V)$ and
$\kappa_\varepsilon=\pi_{x|_{[i]}}(Y\cup Z)$,
where $\kappa_\varepsilon$ does not depend on $V$. Because $C^\varepsilon$ is nonempty, some such set has positive weight, so $\kappa_\varepsilon>0$ and $\pi_y\bigl(\binom{Q_i}{j-q-\varepsilon}\bigr)>0$. The members of $C^\varepsilon$ are exactly the sets $Y\cup Z\cup V$ with $\pi_y(V)>0$. Dividing the last display by its sum over $V$ shows that
\[
\frac{\nu(Y\cup Z\cup V)}{\nu(C^\varepsilon)}=\nu_{y,\,j-q-\varepsilon}(V)
\qquad\Bigl(V\in\tbinom{Q_i}{j-q-\varepsilon}\Bigr),
\]
where both sides vanish when $\pi_y(V)=0$. Therefore
$$
u_1=\sum_{V}\nu_{y,\,j-q-1}(V)\,h(Y\cup\{i\}\cup V) \ \ \text{and } \
u_0=\sum_{W}\nu_{y,\,j-q}(W)\,h(Y\cup W),
$$
the sums running over $\binom{Q_i}{j-q-1}$ and $\binom{Q_i}{j-q}$.

Both layers $j-q-1$ and $j-q$ of $Q_i$ have positive $y$-weight, so Proposition~\ref{prop:transport} gives a nested transport $t$ from $\nu_{y,j-q-1}$ to $\nu_{y,j-q}$ on $Q_i$. Expanding $u_1$ with the first identity in~\eqref{eq:transport} and $u_0$ with the second,
\[
u_1-u_0=\sum_{V\subseteq W}t(V,W)\,\bigl(h(Y\cup\{i\}\cup V)-h(Y\cup W)\bigr).
\]
For each pair $V\subseteq W$, the sets $Y\cup\{i\}\cup V$ and $Y\cup W$ differ by one exchange: $i\notin Y\cup W$, and $W\setminus V$ has exactly one element, hence
$(Y\cup\{i\}\cup V)\setminus(Y\cup W)=\{i\}$ and $(Y\cup W)\setminus(Y\cup\{i\}\cup V)=W\setminus V$.
Hence each bracket has absolute value at most $1$. The numbers $t(V,W)$ are nonnegative and sum to $\sum_V\nu_{y,j-q-1}(V)=1$, so $|u_1-u_0|\le1$.

\emph{Step 2: exponential sums.} For real $\lambda$ and $0\le i\le N$, define
$$
\Phi_i(\lambda)=\sum_{U\in\calS}\nu(U)\,e^{\lambda(M_i(U)-\bar h)}.$$
Then $\Phi_0(\lambda)=1$. We claim that $\Phi_i(\lambda)\le e^{\lambda^2/8}\,\Phi_{i-1}(\lambda)$ for $1\le i\le N$.

Group the sum according to the classes at level $i-1$. A class $C$ as above contributes $\nu(C)e^{\lambda(m-\bar h)}$ to $\Phi_{i-1}(\lambda)$ and
$\nu(C^1)\,e^{\lambda(u_1-\bar h)}+\nu(C^0)\,e^{\lambda(u_0-\bar h)}$
to $\Phi_i(\lambda)$, with a term omitted if its class is empty. If one of $C^1,C^0$ is empty, the other equals $C$, its value of $M_i$ equals $m$ by~\eqref{eq:split-average}, and the two contributions coincide. Otherwise put $p=\nu(C^1)/\nu(C)\in(0,1)$ and $\delta=u_1-u_0$. By~\eqref{eq:split-average}, $m=pu_1+(1-p)u_0$, hence $u_1-m=(1-p)\delta$ and $u_0-m=-p\delta$. The contribution to $\Phi_i(\lambda)$ is therefore
\[
\nu(C)\,e^{\lambda(m-\bar h)}\Bigl(p\,e^{(1-p)\lambda\delta}+(1-p)\,e^{-p\lambda\delta}\Bigr)
\le\nu(C)\,e^{\lambda(m-\bar h)}\,e^{\lambda^2\delta^2/8}
\le\nu(C)\,e^{\lambda(m-\bar h)}\,e^{\lambda^2/8},
\]
by Lemma~\ref{lem:two-point} with $s=\lambda\delta$ and by $|\delta|\le1$ from Step~1. Summing over all classes proves the claim.

\emph{Step 3: the tail bounds.} Iterating Step~2 from $i=N$ down to $i=1$ gives
\begin{equation}\label{eq:mgf}
\sum_{U\in\calS}\nu(U)\,e^{\lambda(h(U)-\bar h)}=\Phi_N(\lambda)\le e^{N\lambda^2/8}\qquad(\lambda\in\R).
\end{equation}
Fix $t>0$ and $\lambda>0$. Every term of the sum in~\eqref{eq:mgf} is nonnegative, and each set $U\in\calS$ with $h(U)\ge\bar h+t$ contributes at least $\nu(U)e^{\lambda t}$. Sets outside $\calS$ have share zero. Hence
\[
\nu\bigl(\{U:h(U)\ge\bar h+t\}\bigr)\le e^{-\lambda t}\,\Phi_N(\lambda)\le\exp\Bigl(-\lambda t+\frac{N\lambda^2}8\Bigr).
\]
Choosing $\lambda=4t/N$ makes the exponent $-4t^2/N+2t^2/N=-2t^2/N$, which proves~\eqref{eq:upper}. The function $-h$ is also $1$-Lipschitz, its weighted average is $-\bar h$, and $\{h\le\bar h-t\}=\{-h\ge-\bar h+t\}$; applying~\eqref{eq:upper} to $-h$ proves~\eqref{eq:lower}. When $t=0$, both bounds state that a share is at most $1$.
\end{proof}

The proof never required $0<x_i<1$: forced and forbidden elements only make some classes have a single branch, and in that case Step~2 contributes the factor $1\le e^{\lambda^2/8}$.

\section*{AI usage declaration}
We used the \texttt{gpt-6-astra} model through xean~\cite{xean},
a harness for AI-assisted mathematical research, to find the counterexample and its proof. Claude Opus 5.5 was used to rewrite the article. All mathematical content has been independently verified by the authors, who take full responsibility for all results.

\end{document}